\documentclass{amsart}

\usepackage{amsmath,amsthm,amssymb}
\usepackage{graphicx}
\usepackage[round]{natbib}

\title{Davis' Convexity Theorem and Extremal Ellipsoids}
\author{Matthias J. Weber \and Hans-Peter Schr\"ocker}
\address{Matthias J. Weber, Hans-Peter Schr\"ocker\\
  Unit Geometry and CAD, University Innsbruck\\
  Technikerstra\ss e 13\\
  A6020 Innsbruck, Austria}
\email{\{matthias.weber, hans-peter.schroecker\}@uibk.ac.at}
\urladdr{http://geometrie.uibk.ac.at}
\keywords{Minimal ellipsoid, maximal ellipsoid, Davis' convexity
  theorem}
\subjclass[2000]
  {52A27, 
   52A20  
  }

\newcommand*{\R}{\mathbb{R}}

\newcommand*{\Norm}[1]{\Vert #1 \Vert}

\DeclareMathOperator{\conv}{conv}

\newtheorem{theorem}{Theorem}
\newtheorem{lemma}[theorem]{Lemma}
\newtheorem{proposition}[theorem]{Proposition}
\theoremstyle{definition}
\newtheorem{definition}[theorem]{Definition}
\theoremstyle{remark}
\newtheorem{remark}[theorem]{Remark}
\begin{document}

\begin{abstract}
  We give a variety of uniqueness results for minimal ellipsoids
  circumscribing and maximal ellipsoids inscribed into a convex
  body. Uniqueness follows from a convexity or concavity criterion on
  the function used to measure the size of the ellipsoid. Simple
  examples with non-unique minimal or maximal ellipsoids conclude this
  article.
\end{abstract}

\maketitle

\section{Introduction}
\label{sec:introduction}

By a classic result in convex geometry the minimal volume ellipsoid
enclosing a convex body $F \subset \R^d$ and the maximal volume
ellipsoid inscribed into $F$ are unique
\citep{john48:_studies_and_essays,danzer57:_loewner_ellipsoid}. Both
ellipsoids are important objects in convex geometry and have numerous
applications in diverse fields of applied and pure mathematics (see
for example
\cite{gruber76:_kennzeichnungen_ellipsoiden,berger90:_convexity} or
the introductory sections of
\cite{kumar05:_minim_volume_ellipsoids,todd07:_khachiyans_algorithm}). More
information on the role of ellipsoids in convex geometry can be found
in \cite{petty83:_convexity_applications} and
\cite[Section~3]{heil93:_special_convex_bodies}.

In this article we are concerned with uniqueness results for minimal
and maximal ellipsoids with respect to size functions different from
the volume. The earliest contribution to this topic is
\cite{firey64:_means_of_convex_bodies} who proved uniqueness of the
minimal quermass integral ellipsoid among all enclosing ellipsoids
with prescribed center. This result can also be deduced from more
general findings of \cite{gruber08:_john_type} and
\cite{schroecker08:_uniqueness_results_ellipsoids}. As to maximal
inscribed ellipsoids we are only aware of
\cite{klartag04:_john_type_ellipsoids} who shows uniqueness with
respect to a vast class of size functions that are defined with the
help of an arbitrary convex body.

In this article we provide uniqueness results for minimal enclosing
and maximal inscribed ellipsoids for further families of size
functions. The basic ideas are similar to that of
\cite{danzer57:_loewner_ellipsoid} and
\cite{schroecker08:_uniqueness_results_ellipsoids}. The new results
are found by applying them to diverse representations of ellipsoids
with the help of symmetric matrices.

After recalling some basic concepts in Section~\ref{sec:preliminaries}
we define the notion of a ``size function'' and, in
Section~\ref{sec:uniqueness-results}, present several different
uniqueness results. In any case it is necessary to study a particular
representation of ellipsoids and properties of an ``in-between
ellipsoid'' in this representation.  Finally, in
Section~\ref{sec:non-unique} we describe a few examples of convex
bodies and size functions with non-unique extremal ellipsoids.

\section{Preliminaries}
\label{sec:preliminaries}

With the help of a positive semi-definite symmetric matrix $A \in
\R^{d \times d}$ and a vector $m \in \R^d$ an ellipsoid can be
described as
\begin{equation}
  \label{eq:1}
  E = \{ x \in \R^d \colon (x-m)^T \cdot A \cdot (x-m) - 1 \le 0 \}.
\end{equation}
The interior of $E$ is the set of all points $x$ that strictly fulfill
the defining inequality. In this article we generally admit degenerate
ellipsoids with empty interior since they may appear as maximal
inscribed ellipsoids. The interior is empty if $A$ is only positive
semi-definite and not positive definite. We call the ellipsoid
\emph{singular} if this is the case and \emph{regular} otherwise.

The vector $m$ is the coordinate vector of the ellipsoid center. A
straight line incident with $m$ and in direction of an eigenvector of
$A$ is called an ellipsoid axis, its semi-axis length $a_i$ is related
to the corresponding eigenvalue $\nu_i$ via $a_i = \nu_i^{-1/2}$.
Since $A$ is symmetric and positive definite, there exist $d$ pairwise
orthogonal axes with real semi-axis lengths.

Note that Equation~\eqref{eq:1} is not the only possibility for
describing ellipsoids. In Section~\ref{sec:uniqueness-results} we will
encounter several alternatives but all of them use a symmetric matrix
and a vector as describing parameters.

There exist different notions for the ``size'' of an ellipsoid. A
natural measure for the size is the ellipsoid's volume, but we may
also take the surface area, a quermass integral, a norm on the vector
of semi-axis lengths etc. More generally, we consider a non-negative
function $f$ on the ordered vector of semi-axis lengths that satisfies
a few basic requirements. By $\R_>$ we denote the set of positive, by
$\R_\ge$ the set of non-negative reals; $\R^d_\ge$ is the set of
vectors $x = (x_1,\ldots,x_d)^T \in \R^d$ with entries $x_i \in
\R_\ge$.

\begin{definition}
  A function $f\colon \R^d_\ge \to \R_\geq$ is called \emph{size function
    for an ellipsoid} if it is continuous, strictly monotone
  increasing in any of its arguments and symmetric, that is, $f(y) =
  f(x)$ whenever $y$ is a permutation of~$x$.
\end{definition}

Denote by $e(A)$ the vector of eigenvalues of a symmetric matrix $A$,
arranged in ascending order. Clearly, $f$ can be extended to the space
of symmetric, positive semi-definite matrices by letting $f(A) = f
\circ e(A)$. Sometimes we will even write $f(E)$ when an ellipsoid $E$
is described by a symmetric matrix~$A$.

Note that $f$ depends only on the eigenvalues of the symmetric matrix
$A$. Hence, it is independent of the position and orientation of $E$.

\section{Uniqueness results}
\label{sec:uniqueness-results}

The uniqueness proofs in this article all follow a certain scheme. We
want to prove that there exists only one minimal enclosing ellipsoid
(with respect to a certain size function $f$) of a convex body $F
\subset \R^d$. Assuming existence of two minimizers $E_0$ and $E_1$ we
construct an ``in-between ellipsoid'' $E_\lambda$ that contains the
common interior of $E_0$ and $E_1$ (and hence also the set $F$) and is
strictly smaller (measured by the size function $f$) than $E_0$ and
$E_1$. Uniqueness results for maximal inscribed ellipsoids can be
obtained in similar fashion.

These type of proof requires the construction of an in-between
ellipsoid $E_\lambda$ that contains $F$ (or is contained in $F$) and
is strictly smaller (or larger) than $E_0$ and $E_1$. Different
constructions of $E_\lambda$ yield different uniqueness results.

\subsection{Image of the unit sphere}
\label{sec:image}

An ellipsoid may be viewed as affine image of the unit ball:
\begin{equation}
  \label{eq:2}
  E = \{ y \in \R^d \colon y = P \cdot x + t,\ x \in \R^d,\ \Norm{x} \leq 1 \},
\end{equation}
where $P \in \R^{d \times d}$ is a (not necessarily regular) matrix
and $t \in \R^d$.

The matrix $P$ is not uniquely determined by the ellipsoid. It is
still possible to apply an automorphic transformation to the unit
sphere before the map $x \mapsto P \cdot x + t$ or an automorphic
transformation to the resulting ellipsoid afterwards. By the left
polar decomposition there exists a symmetric positive semi-definite
matrix $S$ and an orthogonal matrix $U$ such that $P = S \cdot
U$. Hence, we may choose $P$ to be symmetric and positive
semi-definite.

The ordered vector of semi-axis lengths of $E$ is
\begin{equation}
  \label{eq:3}
  a = (a_1, \ldots, a_d)^T = (\nu_1, \ldots, \nu_d)^T,
\end{equation}
where $\nu_i$, $i=1,\dots,d$ are the eigenvalues of $P$. In other
words, we have $a = e(P)$. For reasons that will become clear in the
course of this text we can also write this with the help of the
function
\begin{equation}
  \label{eq:4}
  w^p\colon \R^d \to \R^d,
  \quad
  (x_1,\ldots,x_d)^T \mapsto (\vert x_1 \vert^p,\ldots,\vert x_d \vert^p)^T
\end{equation}
as
\begin{equation}
  \label{eq:5}
  a = w^1 \circ e(P) = e(P).
\end{equation}

\begin{definition}[in-between ellipsoid]
  We define the in-between ellipsoid $E_\lambda$ to two ellipsoids
  $E_0$ and $E_1$ with respect to the representation \eqref{eq:2} as
  \begin{equation}
    \label{eq:6}
    E_\lambda = \{y \in \R^d \colon y = P_\lambda \cdot x + t_\lambda
    ,~ \|x\| \leq 1\},
    \quad
    \lambda \in [0,1]
  \end{equation}
  where
  \begin{equation}
    \label{eq:7}
    E_0 = \{P_0 \cdot x + t_0 \colon \|x\| \leq 1\},
    \quad
    E_1 = \{P_1 \cdot x + t_1 \colon \|x\| \leq 1\},
  \end{equation}
  and
  \begin{equation}
    \label{eq:8}
    P_\lambda = (1-\lambda) P_0 + \lambda P_1,
    \quad
    t_\lambda = (1-\lambda) t_0 + \lambda t_1.
  \end{equation}
\end{definition}

Note that $P_\lambda$ is a symmetric, positive semi-definite matrix
and $E_\lambda$ is indeed an ellipsoid.

\begin{lemma}
  \label{lem:1}
  The in-between ellipsoid $E_\lambda$, $0 \leq \lambda \leq 1$, of
  two ellipsoids $E_0$ and $E_1$ is a subset of the convex hull of the
  two ellipsoids $E_0$ and $E_1$, that is 
  \begin{equation}
    \label{eq:9}
    E_\lambda \subset \conv(E_0, E_1).
  \end{equation}
\end{lemma}

\begin{proof}
  Let $x$ be an element of $E_\lambda$. There exists $y$ with
  $\Norm{y} \leq 1$ such that $x = P_\lambda \cdot y + t_\lambda$. By
  the definition of $P_\lambda$ and $t_\lambda$ we can write
  \begin{equation}
    \label{eq:10}
    x = (1-\lambda) (P_0 \cdot y + t_0) + \lambda (P_1 \cdot y + t_1) =
    (1-\lambda) x_0 + \lambda x_1,
  \end{equation}
  with $x_0 \in E_0$ and $x_1 \in E_1$. Hence, $x$ is in the convex
  hull of $E_0$ and $E_1$ and we conclude $E_\lambda \subset \conv(E_0,
  E_1)$.
\end{proof}

This lemma together with the following proposition already yields a
first uniqueness result for minimal enclosing ellipsoids.

\begin{proposition}[Davis' Convexity Theorem]
  \label{prop:1}
  A convex, lower semi-continuous and symmetric function $f$ of the
  eigenvalues of a symmetric matrix is (essentially strictly) convex
  on the set of symmetric matrices if and only if its restriction to
  the set of diagonal matrices is (essentially strictly) convex.
\end{proposition}

This proposition was stated and proved by
\cite{davis57:_convex_functions} and extended to ``essentially strict
convexity'' by \cite{lewis96:_convex_analysis}.  In
Proposition~\ref{prop:1} ``symmetric'' means that the function $f$ is
independent of the order of its arguments. The precise definition of
``essentially strict convexity'' is rather technical and will be
omitted since we will use only a weaker version of Lewis'
generalization.

We will apply Davis' Convexity Theorem to size functions of
ellipsoids. When proving uniqueness results for minimal ellipsoids we
demand strict convexity of $f$ on $\R_>^d$. For maximal inscribed
ellipsoids we demand strict concavity on $\R_\ge^d$. Results of
\cite{lewis96:_convex_analysis} then guarantee strict
convexity/concavity of $f \circ e$ on the spaces of symmetric matrices
with eigenvalues in $\R_>$, $\R_\ge$, respectively.

\begin{theorem}
  \label{th:1}
  Let $f$ be a size function for ellipsoids such that $f \circ w^1$ is
  strictly concave on $\R_\geq^d$. Further let $F \subset \R^d$ be a
  compact convex body. Among all ellipsoids that are contained in $F$
  there exists a unique ellipsoid that is maximal with respect to~$f$.
\end{theorem}

\begin{proof}
  The existence of a maximal (with respect to $f$) inscribed ellipsoid
  follows from the compactness of $F$ and the continuity of $f \circ
  w^1$. This is explained in great detail in
  \cite{danzer57:_loewner_ellipsoid}.

  To proof uniqueness, we assume existence of two $f$-maximal
  ellipsoids $E_0$ and $E_1$, that is $f(E_0) = f(E_1)$, both
  contained in $F$. We compute the in-between ellipsoid $E_\lambda$
  for $0 < \lambda < 1$ as in \eqref{eq:6}. By Lemma~\ref{lem:1} it is
  contained in the convex hull of $E_0$ and $E_1$ and therefore also
  in $F$. Looking at the size of the in-between ellipsoid we find
  \begin{equation}
    \label{eq:11}
    f(E_\lambda) = f \circ w^1 \circ e(P_\lambda) = f \circ w^1 \circ
    e\big( (1-\lambda) P_0 + \lambda P_1 \big).
  \end{equation}
  Because $P_\lambda$ is a symmetric matrix, we can use Davis'
  Convexity Theorem and find, by strict concavity of $f \circ w^1
  \circ e$,
  \begin{equation}
    \label{eq:12}
    \begin{gathered}
      f \circ w^1 \circ e\big( (1-\lambda) P_0 + \lambda P_1 \big) >
      (1-\lambda) f \circ w^1 \circ e(P_0) + \lambda f \circ w^1 \circ e(P_1) \\
      = f(E_0) = f(E_1).
    \end{gathered}
  \end{equation}
  which is a contradiction.
\end{proof}

\begin{remark}
  The maximal ellipsoids with respect to size functions can be
  computed by a convex program, similar to that described in
  \cite[Section~8.4.2]{boyd04:_convex_optimization}.
\end{remark}

\subsection{Inverse image of the unit sphere}
\label{sec:inverse-image}

In this section we view an ellipsoid as the set
\begin{equation}
  \label{eq:13}
  E = \{ x \in \R^d \colon \|P \cdot x + t\| \leq 1 \},
\end{equation}
where $P \in \R^{d \times d}$ and $t \in \R^d$, that is, as affine
pre-image of the unit ball. Again, it is no loss of generality to
assume that $P$ is symmetric and positive semi-definite. Since we will
use the representation \eqref{eq:13} only for deriving uniqueness
results for minimal ellipsoids we can even assume that $P$ is positive
definite. The ordered vector of semi-axis lengths of $E$ is
\begin{equation}
  \label{eq:14}
  a = w^{-1} \circ e(P).
\end{equation}

\begin{definition}[in-between ellipsoid]
  The in-between ellipsoid $E_\lambda$ to two ellipsoids $E_0$ and
  $E_1$ with respect to the representation \eqref{eq:13} is defined as
  \begin{equation}
    \label{eq:15}
    E_\lambda = \{x \in \R^d \colon \|P_\lambda \cdot x + t_\lambda\| \leq 1\},
    \quad
    \lambda \in [0, 1]
  \end{equation}
  where
  \begin{equation}
    \label{eq:16}
    E_0 = \{x \in \R^d \colon \|P_0 \cdot x + t_0\| \leq 1\},
    \quad
    E_1 = \{x \in \R^d \colon \|P_1 \cdot x + t_1\| \leq 1\},
  \end{equation}
  and
  \begin{equation}
    \label{eq:17}
    P_\lambda = (1-\lambda) P_0 + \lambda P_1,
    \quad
    t_\lambda = (1-\lambda) t_0 + \lambda t_1.
  \end{equation}
\end{definition}
Again, $P_\lambda$ is symmetric and positive definite and $E_\lambda$
is a non-degenerate ellipsoid.

\begin{lemma}
  \label{lem:2}
  Let $E_\lambda$, $0 \leq \lambda \leq 1$, be the in-between
  ellipsoid of two ellipsoids $E_0$ and $E_1$ defined as in
  Equations~\eqref{eq:15}--\eqref{eq:17}. Then the in-between
  ellipsoid $E_\lambda$ encloses the intersection of $E_0$ and~$E_1$.
\end{lemma}

\begin{proof}
  If the intersection of $E_0$ and $E_1$ is empty, nothing has to be
  shown. (Note that this case is irrelevant for the proof of the main
  Theorem~\ref{th:2} below.) Assume therefore that there exists $x \in
  E_0 \cap E_1$, that is,
  \begin{equation}
    \label{eq:18}
    \Norm{P_i \cdot x + t_i} \leq 1,
    \quad
  i \in \{0, 1\}.
  \end{equation}
  We then have
  \begin{equation}
    \label{eq:19}
    1 = (1-\lambda) \cdot 1 + \lambda \cdot 1
    \geq (1-\lambda) \Norm{P_0 \cdot x + t_0} +
    \lambda \Norm{P_1 \cdot x + t_1}.
  \end{equation}
  The triangle inequality implies
  \begin{equation}
    \label{eq:20}
    \begin{aligned}
    & (1-\lambda) \Norm{P_0 \cdot x + t_0} + \lambda \Norm{P_1 \cdot x
      + t_1 } \geq \\
    & \Norm{(1-\lambda) \bigl( P_0 \cdot x + t_0 \bigr) + \lambda
      \bigl( P_1 \cdot x + t_1 \bigr) } =\\
    & \Norm{ \big( (1-\lambda) P_0 + \lambda P_1 \big) \cdot x + \big(
    (1-\lambda) t_0 + \lambda t_1 \big) } =\\
    & \Norm{ P_\lambda \cdot x + t_\lambda }.
    \end{aligned}
  \end{equation}
  Combining \eqref{eq:19} and \eqref{eq:20} we see that $\Norm{
    P_\lambda \cdot x + t_\lambda } \leq 1$. This shows that $x \in
  E_\lambda$. Hence $E_0 \cap E_1 \subset E_\lambda$ and the proof is
  complete.
\end{proof}

\begin{theorem}
  \label{th:2}
  Let $f$ be a size function for ellipsoids such that $f \circ w^{-1}$
  is strictly convex on $\R_>^d$. Further let $F \subset \R^d$ be a
  compact convex body. Among all ellipsoids that contain $F$ there
  exists a unique ellipsoid that is minimal with respect to~$f$.
\end{theorem}

\begin{proof}
  The existence of a minimal (with respect to $f$) ellipsoid that
  encloses $F$, follows from the compactness of $F$ and the continuity
  of $f \circ w^{-1}$ (see again \cite{danzer57:_loewner_ellipsoid}).

  To proof uniqueness, we assume existence of two $f$-minimal
  ellipsoids $E_0$ and $E_1$, that is $f(E_0) = f(E_1)$, both
  containing $F$. We compute the in-between ellipsoids $E_\lambda$ for
  $0 < \lambda < 1$, as in \eqref{eq:15}. By Lemma~\ref{lem:2} it
  contains the common interior of $E_0 \cap E_1$ and hence also
  $F$. Looking at the size of $E_\lambda$ we find
  \begin{equation}
    \label{eq:21}
    f(E_\lambda) = f \circ w^{-1} \circ e(P_\lambda) = f
    \circ w^{-1} \circ e\big( (1-\lambda) P_0 + \lambda P_1
    \big).
  \end{equation}
  Because $P_\lambda$ is a symmetric matrix, we can use Davis'
  Convexity Theorem (see Proposition~\ref{prop:1} on
  page~\pageref{prop:1}). It implies that $f \circ w^{-1} \circ e$ is
  strictly convex. Therefore we can write
  \begin{equation}
    \label{eq:22}
    f \circ w^{-1} \circ e\big( (1-\lambda) P_0 + \lambda P_1 \big) <
    (1-\lambda) f \circ w^{-1} \circ e(P_0) + \lambda f \circ w^{-1}
    \circ e(P_1).
  \end{equation}
  Because $E_0$ and $E_1$ have the same size it follows that
  \begin{equation}
    \label{eq:23}
    f(E_\lambda) = f \circ w^{-1} \circ e(P_\lambda) < f \circ w^{-1}
    \circ e(P_0) = f(E_0) = f(E_1).
  \end{equation}
  We have now that the size of $E_\lambda$ is smaller than the size of
  $E_0$ and $E_1$. Together with Lemma~\ref{lem:2} this constitutes a
  contradiction to the assumed minimality of $E_0$ and $E_1$ and
  finishes the proof.
\end{proof}

\subsection{Extremal affine images of convex unit balls}
\label{sec:extremal-affine-images}

It is easy to see that the proves of Theorems~\ref{th:1} and
\ref{th:2} remain true if we replace the Euclidean unit ball by an
arbitrary centrally symmetric convex body, centered at the origin, and
measure its size by the volume. Hence, we can state a much more
general result:

\begin{theorem}
  The volume-minimal circumscribing affine image of an arbitrary convex
  unit ball to a compact convex body $F$ is unique. The same is true
  for volume-maximal inscribed affine image of an arbitrary convex
  unit ball.
\end{theorem}

\subsection{Algebraic equation}
\label{sec:equation}

The maybe most straightforward way to represent a non-degenerate
ellipsoid $E \subset \R^d$ uses the algebraic equation of $E$:
\begin{equation}
  \label{eq:24}
  E  = \{x \in \R^d \colon (x-m)^T \cdot A \cdot (x-m) \leq 1 \},
\end{equation}
with a symmetric, positive definite matrix $A \in \R^{d \times d}$ and
$m \in \R^d$.  The vector $a$ of ordered semi-axis lengths of $E$ is
found as
\begin{equation}
  \label{eq:25}
  a = w^{-1/2} \circ e(A).
\end{equation}

The equation of $E$ can also be written with the help of a single
matrix of dimension $(d+1) \times (d+1)$:
\begin{equation}
  \label{eq:26}
  E = \{X \in \R^{d+1} \colon X^T \cdot M \cdot X \leq 0\},
\end{equation}
where
\begin{equation}
  \label{eq:27}
  X =
  \begin{pmatrix}
    1\\
    x
  \end{pmatrix},
  \quad
  M =
  \begin{pmatrix}
    -1 & -m^T \cdot A'\\
    -A' \cdot m & A'
  \end{pmatrix}
  \quad \text{and} \quad
  A' = \frac{A}{1-m^T \cdot A \cdot m}.
\end{equation}
If we define the in-between ellipsoid $E_\lambda$ to two ellipsoids
$E_0$ and $E_1$ with respect to the representation \eqref{eq:26} by
building a convex sum of the two homogeneous matrices that define
$E_0$ and $E_1$,
\begin{equation}
  \label{eq:29}
  E_\lambda = \{X \in \R^{d+1} \colon X^T \cdot M_\lambda \cdot X \leq 0 \},
  \quad
  \lambda \in [0,1]
\end{equation}
where
\begin{equation}
  \label{eq:30}
  M_\lambda = (1-\lambda) M_0 + \lambda M_1,
\end{equation}
we arrive at the situation discussed in
\cite{schroecker08:_uniqueness_results_ellipsoids}. The main
uniqueness result is

\begin{proposition}
  \label{th:3}
  Let $f$ be a size function and $f \circ w^{-1/2}$ be a strictly
  convex function on $\R_>^d$. Further let $F \subset \R^d$ be a
  compact convex body. Among all ellipsoids that contain $F$ there
  exists a unique ellipsoid that is minimal with respect to~$f$.
\end{proposition}

\subsection{Dual equation}
\label{sec:dual-equation}

An ellipsoid can also be viewed as the set of hyperplanes that
intersect the (point-set) ellipsoid in real points. Using hyperplane
coordinates, this description is formally the same as in
Section~\ref{sec:equation}:
\begin{equation}
  \label{eq:32}
  E = \{ u \in \R^d \colon (u-c)^T \cdot B \cdot (u-c) \leq 1 \},
\end{equation}
where $B \in \R^{d \times d}$ is a symmetric, positive semi-definite
matrix and $c \in \R^d$. In homogeneous form this is
\begin{equation}
  \label{eq:33}
  E = \{ U \in \R^{d+1} \colon U^T \cdot N \cdot U \leq 0 \},
\end{equation}
where
\begin{equation}
  \label{eq:34}
  U =
  \begin{pmatrix}
    1\\
    u
  \end{pmatrix},
  \quad
  N =
  \begin{pmatrix}
    -1 & -c^T \cdot B'\\
    -B' \cdot c & B'
  \end{pmatrix},
  \quad \text{and} \quad
  B' = \frac{B}{1 - c^T \cdot B \cdot c}.
\end{equation}
Translating the center of $E$ to the origin, this description becomes
\begin{equation}
  \label{eq:35}
  E_o = \{ U \in \R^{d+1} \colon U^T \cdot N_o \cdot U \leq 0 \},
\end{equation}
where
\begin{equation}
  \label{eq:36}
  N_o =
  \begin{pmatrix}
    -1 & 0^T\\
    0 & B' \cdot c \cdot c^T \cdot B' + B'
  \end{pmatrix}.
\end{equation}
In this representation, the vector of semi-axis lengths is
\begin{equation}
  \label{eq:37}
  a = w^{1/2} \circ e(B' \cdot c \cdot c^T \cdot B' + B').
\end{equation}

\begin{definition}[in-between ellipsoid]
  We define the in-between ellipsoid $E_\lambda$ to two ellipsoids
  $E_0$ and $E_1$ with respect to the representation \eqref{eq:33} by
  building the convex sum of the two defining homogeneous matrices:
  \begin{equation}
    \label{eq:38}
    E_\lambda = \{ U \in \R^{d+1} \colon U^T \cdot N_\lambda \cdot U \leq 0  \},
    \quad
    \lambda \in [0,1]
  \end{equation}
  where
  \begin{equation}
    \label{eq:39}
    N_\lambda = (1-\lambda) N_0 + \lambda N_1.
  \end{equation}
\end{definition}

Note that we have no guarantee that $E_\lambda$ is really an ellipsoid
for all values $\lambda \in [0,1]$. It is, however, an ellipsoid at
least in the vicinity of $\lambda = 0$ and $\lambda = 1$ and this is
all we need. For reasons of simplicity we will not always mention this
explicitly and still refer to $E_\lambda$ as ``in-between ellipsoid''.

\begin{lemma}
  \label{lem:4}
  The in-between ellipsoid $E_\lambda$ of two ellipsoids $E_0$ and
  $E_1$ lies inside the convex hull of $E_0$ and $E_1$, that is
  \begin{equation}
    \label{eq:40}
    E_\lambda \subset \conv(E_0, E_1),
  \end{equation}
  at least for values of $\lambda$ in the vicinity of $0$ and~$1$.
\end{lemma}

In order to prove Lemma~\ref{lem:4} it is sufficient to consider the
case $d=2$. This can be seen as follows: Let $x$ be a point in
$E_\lambda$ and take a plane $\pi$ through $x$ and the centers of
$E_0$ and $E_1$, respectively. The in-between ellipsoid $E_\lambda$
intersects $\pi$ in an ellipse $E'_\lambda$ that is obtained as
in-between ellipse to $\pi \cap E_0$ and $\pi \cap E_1$. Hence, $x$
lies in $E_\lambda$ if and only if it lies in~$E'_\lambda$.

The proof for $d=2$ can be carried out by straightforward
computation. It requires, however, a case distinction, is rather
technical and does not provide useful insight. Therefore, we omit it
at this place. It will be published in the first author's doctoral
thesis.

\begin{theorem}
  \label{th:5}
  Let $f$ be a size function for ellipsoids such that $f\circ w^{1/2}$
  is strictly concave on $\R_\geq^d$. Further let $F \subset \R^d$ be
  a compact convex body. Among all ellipsoids with a fixed center that
  are inscribed into $F$ there exists a unique ellipsoid that is
  maximal with respect to~$f$.
\end{theorem}

Once we have realized that we can describe $E_0$ and $E_1$ by
homogeneous matrices
\begin{equation}
  \label{eq:41}
  N_i =
  \begin{pmatrix}
    -1 & 0^T\\
    0 & B_i
  \end{pmatrix},
  \quad i=0,1
\end{equation}
the proof is quite similar to the proof of Theorem
\ref{th:1}.

\begin{remark}
  The uniqueness results of Theorems~\ref{th:1}, \ref{th:2}, and
  \ref{th:5} also hold if we look for extremal ellipsoids only among
  ellipsoids with prescribed axes. Theorems~\ref{th:1} and \ref{th:2}
  remain true if the center is prescribed.
\end{remark}

\section{Non-uniqueness results}
\label{sec:non-unique}

In this section we give two simple examples of size functions and
convex sets such that the corresponding extremal ellipsoids are not
unique. In view of our results, the size functions lack a convexity or
concavity property. While non-uniqueness in both examples is rather
obvious we feel the need to publish them since we are not aware of a
single similar counter-example. Only
\cite{behrend38:_kleinste_ellipse} mentions the non-uniqueness of
maximal inscribed circles. A trivial example is two congruent circles
inscribed into their convex hull.

\subsection*{Minimal ellipsoids with non-convex size function}
\label{sec:examples_counter-minell}

\begin{figure}
  \centering
  \begin{minipage}[t]{0.5\linewidth}  
    \rule{\linewidth}{0pt}
    \centering
    \includegraphics{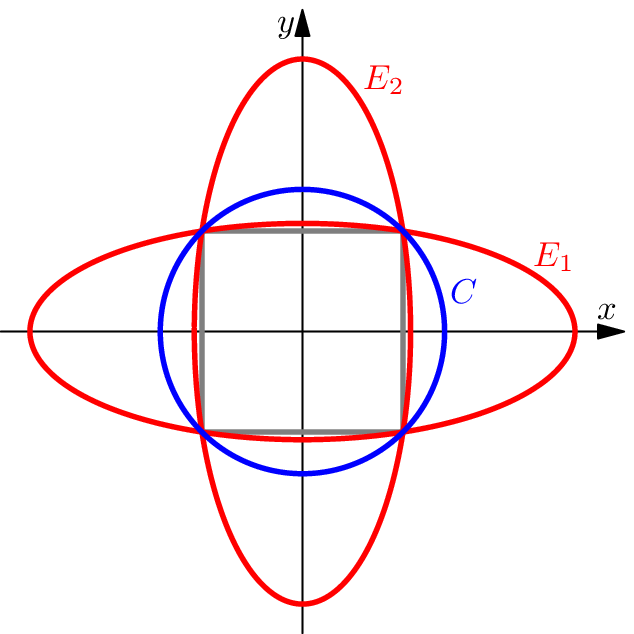}
  \end{minipage}%
  \begin{minipage}[t]{0.5\linewidth}  
    \rule{\linewidth}{0pt}
    \centering
    \includegraphics{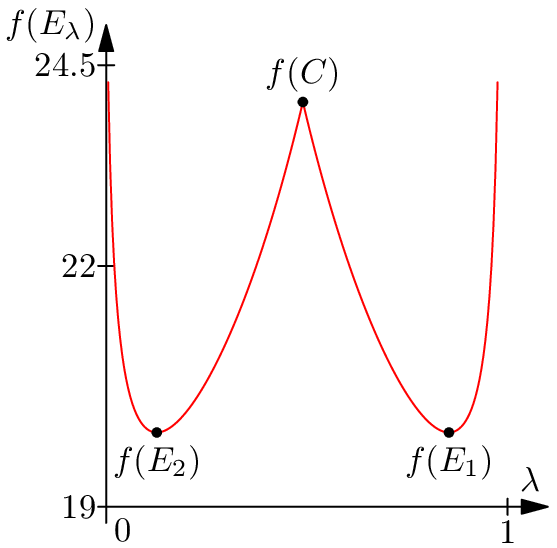}
  \end{minipage}
  \caption{Non-unique minimal ellipsoids through the vertices of a
    square.}
  \label{fig:minEllGegenbsp}
\end{figure}

Denote by $F \subset \R^2$ the set of four points with coordinates
$(\pm 1,\pm 1)$ and let $f$ be the non-convex size function
\begin{equation*}
  f \colon \R_\geq^2 \to \R_\geq, ~(a,b) \mapsto \max\{a,b\} + 16 \min\{a,b\}.
\end{equation*}
If the $f$-minimal ellipse to $F$ was unique it must have four axis
of symmetry and therefore it must be the circle $C$ through the points
of $F$. But the size of the two ellipses $E_1$ and $E_2$
\begin{equation*}
  \begin{aligned}
    E_1\colon& \Bigl(\frac{32}{257} 2^{\frac{1}{3}} - \frac{4}{257}
    2^{\frac{2}{3}} + \frac{1}{257}\Bigr) x^2 + 
    \Bigl(-\frac{32}{257} 2^{\frac{1}{3}} + \frac{4}{257}
    2^{\frac{2}{3}} + \frac{256}{257}\Bigr) y^2 - 1 \le 0\\ 
    E_2\colon& \Bigl(-\frac{32}{257} 2^{\frac{1}{3}} + \frac{4}{257}
    2^{\frac{2}{3}} + \frac{256}{257}\Bigr) x^2 + 
    \Bigl(\frac{32}{257} 2^{\frac{1}{3}} - \frac{4}{257}
    2^{\frac{2}{3}} + \frac{1}{257}\Bigr) y^2 -1 \le 0 
  \end{aligned}
\end{equation*}
is smaller than the size of the circle (compare
Figure~\ref{fig:minEllGegenbsp}):
\begin{equation*}
  f(E_1) = f(E_2) \approx 19.9248 < f(C) \approx 24.0416
\end{equation*}
The ellipses $E_1$ and $E_2$ are the minimizers of $f$ among all
ellipses $E_\lambda$ through the four points of $F$.
Figure~\ref{fig:minEllGegenbsp}, right, displays the plot of the size
function for all ellipses in the pencil of conics spanned by these
points.

\subsection*{Maximal ellipsoids with non-concave size function} 
\label{sec:examples_counter-maxell}

\begin{figure}
  \begin{minipage}[b]{0.5\linewidth}
    \rule{\linewidth}{0pt}
    \centering
    \includegraphics{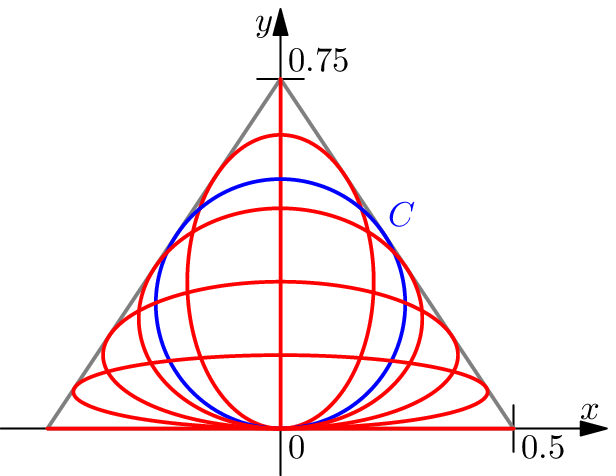}
  \end{minipage}%
  \begin{minipage}[b]{0.5\linewidth}
    \rule{\linewidth}{0pt}
    \centering
    \includegraphics{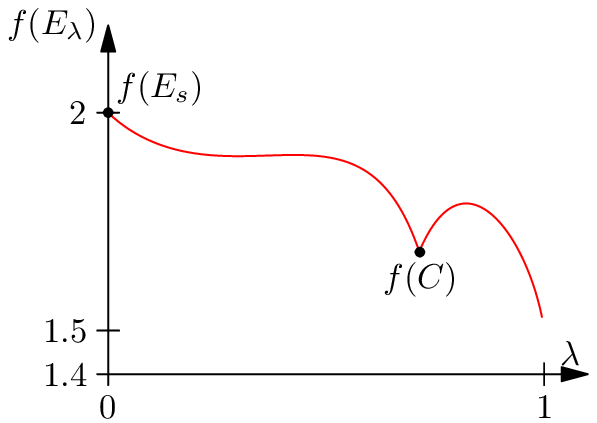}
  \end{minipage}
  \caption{The arc-length of some ellipses inscribed into an
    equilateral triangle}
  \label{fig:maxEllsizefunc}
\end{figure}

Let $F \subset \R^2$ be the equilateral triangle with side length
$1$ (see Figure~\ref{fig:maxEllsizefunc}). The size function under
consideration is the arc-length of an
ellipse.  We will demonstrate that the inscribed ellipse of maximal
arc length is not unique.  This is particularly interesting since the
minimal arc-length enclosing ellipse is known to be unique, see
\cite{firey64:_means_of_convex_bodies,gruber08:_john_type,schroecker08:_uniqueness_results_ellipsoids}.

The arc-length of an ellipse with semi-axis length $a$ and $b$ can be
expressed in terms of the complete elliptic integral of first kind
\begin{equation*}
  f(a, b) = 4 \max\{a,b\} E(1 - \frac{\min\{a,b\}}{\max\{a,b\}})
  \quad\text{where}\quad
  E(k) = \int_0^1 \frac{\sqrt{1-k^2 t^2}}{\sqrt{1-t^2}}\;\mathrm{d}t.
\end{equation*}

If the maximal arc-length ellipse contained in $F$ was unique it must
share the triangle's symmetries. Therefore, it must be the in-circle
$C$. But the arc-length of the ellipsoid $E_s$ that degenerates to the
triangle side on the $x$-axis is greater than that of the circle:
$f(E_s) = 2 > \pi/\sqrt{3} = f(C)$, see
Figure~\ref{fig:maxEllsizefunc}. This shows that the maximal
arc-length ellipse inscribed into an equilateral triangle is not
unique. The plot in Figure~\ref{fig:maxEllsizefunc}, right, depicts
the size function of the drawn inscribed ellipses. The circle
corresponds to the kink in the graph.

\section{Conclusion}
\label{sec:conclusion}

We studied uniqueness results of minimal circumscribed and maximal
inscribed ellipsoids. Uniqueness can be guaranteed if the function
used for measuring the ellipsoid size satisfies a certain convexity or
concavity condition. Summarizing our findings we can state that the
minimal enclosing ellipsoid with respect to a size function $f$ is
unique if $f \circ w^p$ is convex for $p \in \{-1, -1/2\}$. The maximal
inscribed ellipsoid is unique if $f \circ w^p$ is concave for $p = 1$
or for $p=1/2$ if the center is prescribed. Uniqueness for $p=1/2$
under general assumptions is still an open question.

\section*{Acknowledgments}

The authors gratefully acknowledge support of this research by the
Austrian Science Foundation FWF under grant P21032.


\end{document}